\documentclass[psamsfonts]{amsart}
\input epsf.tex
\usepackage{a4wide}
\font\tenbb=msbm10 scaled 1000
\font\sevenbb=msbm7 scaled 1000
\font\fivebb=msbm5 scaled 1000
\newfam\bbfam
\textfont\bbfam=\tenbb
\scriptfont\bbfam=\sevenbb
\scriptscriptfont\bbfam=\fivebb

\input{psfig.tex}
\usepackage{latexsym,amsfonts,amssymb,amsbsy}
\usepackage{amsmath} 

\newtheorem{thm}{Theorem}[section] 
\newtheorem{prop}[thm]{Proposition} 
\newtheorem{cor} [thm]{Corollary} 
\newtheorem{lem}[thm] {Lemma} 

\newtheorem{df} {Definition}
\newtheorem{quest}{Question}

\newcommand{\edim}{\hspace*{\fill}$\Box$ \vspace{10pt}}

\newcommand{\bns}{B_n(\Sigma_{g,p})}

\newcommand{\gend}{\sg_1,\dots,\sg_{n-1},
a_1, \dots,a_{g}, b_1, \dots ,b_g,  z_1, \dots, z_{p-1}} 
\newcommand{\genu}{\sigma_1, \dots, \sigma_{n-1}}
\newcommand{\genh}{\sg_1,\dots,\sg_{n-1},
a_1, \dots, a_{g}, b_1, \dots , b_g,  z_1, \dots, z_{p-1}} 

\newcommand{\sg}{\sigma}

\newcommand{\sjp}{\sg_{j}}
\newcommand{\sip}{\sg_{i}}
\newcommand{\siip}{\sg_{i+1}}

\newcommand{\iji}{\sip \, \siip \, \sip}
\newcommand{\jij}{\siip \, \sip \, \siip}
\newcommand{\ij}{\sip \, \sjp}
\newcommand{\ji}{\sjp \, \sip}

\begin{document}
\title{Questions on  surface braid groups}

\author{Paolo Bellingeri}  
\author{Eddy Godelle}
\address{Univ. Pisa, Dipartimento di Matematica, 56205 Pisa, Italy}  
 \email{bellingeri@mail.dm.pisa.it}
\address{Univ. Caen, Laboratoire de Math\'ematiques Nicolas Oresme, 14032
  Caen, France}  
 \email{Eddy.Godelle@math.unicaen.fr}
%\thanks{\
%\mbox{\hspace{11pt}}{\it Mathematics Subject Classification:}
%Primary: 20F36.\\
%\mbox{\hspace{11pt}}{\it Key words}: 
%Braids, surface braids, group presentations

\begin{abstract}
We provide new group presentations for surface braid groups
which are positive. We study some properties of such presentations
and we solve the conjugacy problem in a particular case.
%In this paper we ask some questions on our group presentations 
% that seem to us relevant% to
%better understand algebraical properties of surface braid groups. 
\end{abstract}
 
\maketitle
%%%%%%%%%%%%%%%%%%%%%%%%%%%%%%%%%%%%%%%%%%%%%%%%%%%%%%%%%%%%%%%%%  

%%%%%%%%%%%%%%%%%%%%%%%%%%%%%%%%%%%%%%%%%%%%%%%%%%%%%%%%%%%%%%%%% 

\section{Introduction and Motivation}

Let $\Sigma_{g,p}$ be an orientable surface of  genus $g$ with $p$
  boundary components.
  For instance, 
 $\Sigma_{0,0}$ is the 2-sphere, $\Sigma_{1,0}$ is the
 torus, 
and $\Sigma_{0,1}$ corresponds to the disk.

 A \emph{geometric braid} on $\Sigma_{g,p}$ based at 
$\mathcal{P}$ is a collection $B=(\psi_1, \dots, \psi_n)$ of $n\geq 2$
 paths from $[0,1]$ to $\Sigma_{g,p}$   such that $\psi_i(0)= P_i$, $\psi_i(1)
 \in \mathcal{P}$
and $\{ \psi_1(t), \dots, \psi_n(t)\}$ are distinct points  for all $t
\in [0, 1]$.
Two braids are considered as equivalent if they are isotopic.
 The usual  product of paths defines a group structure
 on the equivalence classes  of  braids.
 This group doesnot depend, up to isomorphism,  on the choice of $\mathcal{P}$.
 It is called \emph{the (surface) braid group on $n$ strands} on
 $\Sigma_{g,p}$ and denoted by $\bns$.  The group $B_n(\Sigma_{0,1})$ is the classical braid group $B_n$ on $n$ strings.
 Some elements of $\bns$ are shown in Figure \ref{gen:sb1}.
The braid $\sg_i$ corresponds to the 
standard generator of $B_n$ and it can be represented by a geometric braid on
$\Sigma_{g,p}$ where all the strands are trivial except the $i$-th one and 
the $(i+1)$-th one. The  $i$-th strand goes from $P_i$ to $P_{i+1}$ and the $(i+1)$-th strand goes from
$P_{i+1}$ to $P_{i}$ according to Figure 1. 
The loops  $\delta_1, \dots, \delta_{p+2g-1}$ based on $P_1$ in Figure 1 represent standard generators of $\pi_1(\Sigma_{g,p})$. 
By its definition,  $B_1(\Sigma_{g,p})$ is isomorphic to
$\pi_1(\Sigma_{g,p})$. We can also consider 
$\delta_1, \dots, \delta_{p+2g-1}$ as
braids on $n$ strands on $\Sigma$, where last $n-1$ strands are trivial.

It is well known since E. Artin (\cite{Art}) that the braid group
$B_n$ has a positive presentation 
(see  for instance \cite{MurKur} Chapter 2 Theorem 2.2), i.e.  a
group presentation which involves only  generators and not their
inverses. Hence one can associate a (braid) monoid $B_n^+$ with the
same presentation, but as a monoid presentation. It turns out that the
braid monoid $B_n^+$ is a Garside monoid (see \cite{Deh}), that is a
monoid with a good divisibility structure, and that the braid group
$B_n$ is the group of fractions of the monoid $B_n^+$. As a
consequence, the natural morphism of monoids from $B_n^+$ to $B_n$ is
into, and we can solve the word problem, the conjugacy problem and
obtain normal forms in $B_n$(see \cite{BrS,Cha,Deh,Del,Gar}). These
results extend to Artin-Tits groups of spherical type
which are a well-known algebraic generalization of the braid group $B_n$
(\cite{BrS,Cha,Deh,Del}).\\ 
In the case of surface braid groups $\bns$, some group
presentations are known but they are not positive. 
Furthermore, questions as the conjugacy problem are not solved in the
general case. 
The word problem in surface braid groups is known to be solvable (see
\cite{Gon}) 
even if  algorithms are not as efficient as the ones proposed for the
braid group $B_n$.\\ 
In this note we provide positive presentations for $\bns$ and we
address questions related
 to the conjugacy problem of surface braid groups.
 We do not discuss
 the case of $B_n(\Sigma_{0,0})$, the braid group  
on the $2$-sphere; this is a particular case with
specific properties. 
For instance if $\Sigma$ is an oriented surface, the  surface braid
group 
$B_n(\Sigma)$  has
torsion 
elements only and only if $\Sigma$ is the $2$-sphere (see \cite{GilBus} page 277, \cite{FadBus} page 255, and \cite{ParRol} proposition 1.5).\\ 
In Section 2 and 3 we focus on braid groups on surfaces with boundary
 components and without 
 boundary components respectively. In Section 4, we investigate the
 special
 case of $B_2(\Sigma_{1,0})$ and  we solve the word problem
 and the
 conjugacy problem for this group.   

\vspace{5pt}

\medskip 
\epsfysize=7truecm
\epsfxsize=13truecm
\medskip
\centerline{\epsfbox{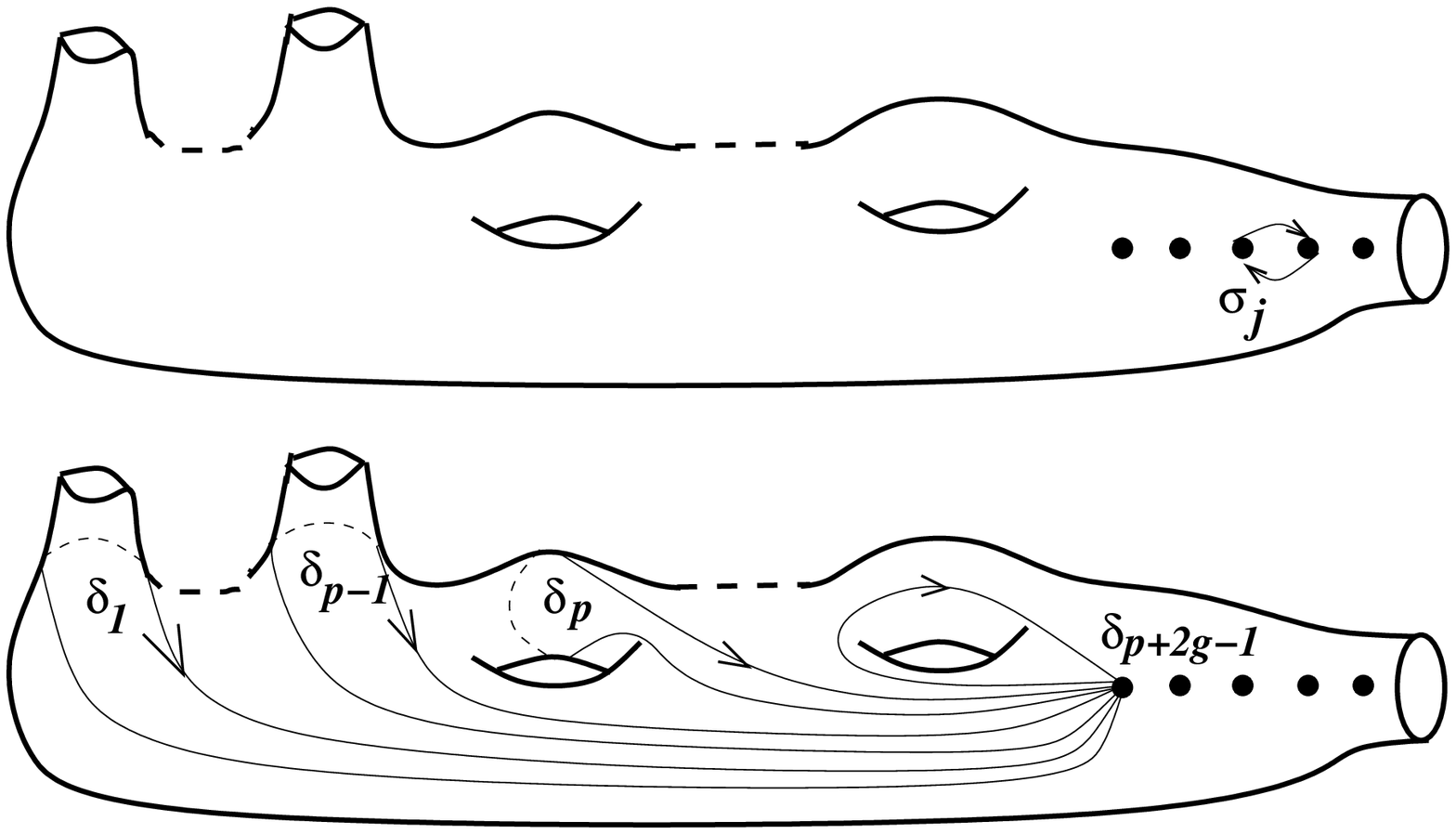}}
\label{gen:sb1}
\medskip
\centerline{Figure 1: some braid elements}
\medskip

\section{Braid groups on surfaces with boundary components}
In this section we investigate  braid groups on oriented surfaces with a 
positive number of boundary components. Our first objective is to prove Theorem \ref{thm:presbg1}:
\begin{thm} \label{thm:presbg1}  Let $n$ and $p$ be positive integers. Let $g$ be a non negative integer. 
Then, the group $\bns$ admits the following group presentation:\\
$\bullet$ Generators: $\genu,\delta_1,\cdots , \delta_{2g+p-1}$;\\
$\bullet$ Relations:\\
-Braid relations:\\
\begin{tabular}{lll}
$(BR1)$ &$\ij  =  \ji$&for $| i-j | \ge 2$;\\
$(BR2)$ &$\iji  = \jij$&for $1 \le i\le n-1$.\end{tabular}\\
-  Commutative relations between surface braids:\\
\begin{tabular}{lll}
$(CR1)$ &$\delta_r \sg_i= \sg_i \delta_r$  &for $i\not= 1$; $1 \le r \le 2g+p-1$;\\
$(CR2)$ &$\delta_{r}  \sigma_1 \delta_{r} \sigma_1 =  \sigma_1 \delta_{r} \sigma_1 \delta_{r}$& $1\leq r\leq 2g+p-1$;\\
$(CR3)$&$\delta_r \sigma_1 \delta_r \delta_{s} \sigma_1 =  \sigma_1 \delta_r\delta_{s} \sigma_1 \delta_r$ & 
for $1 \le r<s \le 2g+p-1$ with\\&& $(r,s)\neq (p+2i, p+2i+1)$, $0\leq i \leq g-1$. \end{tabular}\\
-  Skew commutative relations on the handles:\\
$(SCR1)$  $\sg_1 \delta_{r+1} \sigma_1 \delta_{r} \sigma_1=  \delta_{r} \sigma_1 \delta_{r+1}$  for $r = p+2i$ where $0 \le i \le g-1$.
\end{thm}
The above presentation can be compared to the presentation of $B_{g,n}$ given in \cite{HaeLam} page 18.
\proof
Let us denote by $\widetilde{\bns}$ the group defined by the presentation given in Theorem~\ref{thm:presbg1}. 
We prove that the group $\widetilde{\bns}$ is isomorphic to the  group $\bns$ using the presentation given in 
Theorem~\ref{thm:pres}. Let $\psi:\{\gend\}\to 
\{\genu,\delta_1,\cdots , $ $\delta_{2g+p-1}\}$ be the set-map defined by $\psi(\sg_i)=\sg_i$ 
for $i=1, \dots, n-1$, $\psi(a_r)=\delta_{p+2(r-1)}^{-1}$, $\psi(b_r)=\delta_{p+2(r-1)+1}^{-1}$ 
for $r=1, \dots, g$ and $\psi(z_j)=\delta_j^{-1}$ for  $j=1, \dots, p-1$. We claim that $\psi$ extends to a 
homomorphism of groups $\psi: \bns \to \widetilde{\bns}$. We have to verify that the image by $\psi$ of the 
braid relations and of the relations of type (R1)-(R8) are true in $\widetilde{\bns}$. It is enough to verify that 
\stepcounter{equation}
\begin{eqnarray} \label{eqaprouver}\sigma_1^{-1} \delta_r^{-1} \sigma_1 \delta_s^{-1} = 
\delta_s^{-1} \sigma_1^{-1} \delta_r^{-1} \sigma_1\end{eqnarray}  for $1\leq r<s\leq 2g+p-1$ and $(r,s)\neq (p+2k,p+2k+1)$
 which corresponds to the image by $\psi$ in $\widetilde{\bns}$ of the relations of type  (R3), (R6) and (R7); the other cases 
are true as they are relations of the presentation of $\widetilde{\bns}$.\\
The relations of type (CR3) can be written $\delta_r \sg_1 \delta_r \delta_s  = \sg_1 \delta_r \delta_s \sg_1 \delta_r \sg_{1}^{-1}$. 
From the relations of type (CR1) we deduce that, in $\widetilde{\bns}$, the equalities 
$\delta_r\sg_1 \delta_r\delta_s =\delta_r\sg_1\delta_r\sg_1\sg_{1}^{-1}\delta_s =  \sg_1 \delta_r\sg_1\delta_r\sg_{1}^{-1}\delta_s$ 
holds. Hence we obtain $\sg_1\delta_r\delta_s\sg_1\delta_r\sg_{1}^{-1}=\sg_1\delta_r\sg_1\delta_r  \sg_{1}^{-1}\delta_s$.   
From this equality, we derive that $\delta_s \sg_1 \delta_r \sg_{1}^{-1}=\sg_1 \delta_r  \sg_{1}^{-1} \delta_s$, and finally 
we get the relations (\ref{eqaprouver}).
        
        On the other hand, consider $\overline{\psi}$ the set-map defined from $\{\genu,\delta_1,\cdots , \delta_{2g+p-1}\}$ to 
$\{\genh\}$ by $\overline{\psi}(\sg_i)=\sg_i$ for $i=1, \dots, n-1$,
        $\overline{\psi}(\delta_j)=z_j^{-1}$ for  $j=1, \dots, p-1$, $\overline{\psi}(\delta_{p+2(r-1)})=a_r^{-1}$ and  
        $\overline{\psi}(\delta_{p+2(r-1)+1)}=b_r^{-1}$ for $r=1, \dots, g$.
        We prove that $\overline{\psi}$ extends to an homomorphism of groups from $\widetilde{\bns}$ to $\bns$.
        Since braid relations and  the images by $\overline{\psi}$ of
        the relations of type $(ER)$, $(CR1)$ and 
 $(CR2)$ are verified, it suffices to check that the equalities 
corresponding to relations of type $(CR3)$ hold in $\bns$.
        We verify that the equality $a_r^{-1} \sigma_1 a_r^{-1} a_{s}^{-1} \sigma_1 =  \sigma_1 a_r^{-1} a _{s}^{-1} 
\sigma_1 a_r^{-1}$ for  $(1 \le r>s \le g)$ holds in $\bns$. The other cases can easily be verified by the reader. 
From the relations of type (R2), it follows that $ a_r^{-1} \sigma_1 a_r^{-1} =  \sg_1 a_r^{-1} \sigma_1 a_r^{-1} \sg_1^{-1}$; 
thus we have $a_r^{-1} \sigma_1 a_r^{-1} a_{s}^{-1} \sigma_1 = \sg_1 a_r^{-1} \sigma_1 a_r^{-1} \sg_1^{-1} a_{s}^{-1} \sigma_1$. 
Applying  relations of type (R3) we deduce that $ a_r^{-1} \sg_1^{-1} a_{s}^{-1} \sigma_1= \sg_1^{-1} a_{s}^{-1} \sigma_1   a_r^{-1}$  
and therefore the equalities $a_r^{-1} \sigma_1 a_r^{-1} a_{s}^{-1} \sigma_1 = \sigma_1 a_r^{-1} a _{s}^{-1} \sigma_1 a_r^{-1}$ 
hold in $\bns$.
        Then, the morphism $\overline{\psi}$ from $\widetilde{\bns}$ to $\bns$ is well defined and it is the 
inverse of $\psi$. Hence, $\bns$ is isomorphic to $\widetilde{\bns}$.  
\edim

We remark that the presentation given in Theorem \ref{thm:presbg1} is positive and has less types of relations than the 
presentation given in Theorem \ref{thm:pres}. 
\begin{lem} \label{lem:reduction} Let $G$ be a group and let
  $\sigma,\delta,\delta'$ be in $G$.\\(i) If
 (a) $\delta(\sigma\delta\delta'\sigma) = (\sigma\delta\delta'\sigma)\delta$, 
(b) $\sigma\delta\sigma\delta = \delta\sigma\delta\sigma$ and (c) $\sigma\delta'\sigma\delta' = \delta'\sigma\delta'\sigma$ 
then $\delta'(\sigma\delta\delta'\sigma) = (\sigma\delta\delta'\sigma)\delta'$.\\
(ii) If (a) $\delta(\sigma\delta\delta'\sigma) = \sigma(\delta\delta'\sigma\delta)$, 
(b) $\delta'(\sigma\delta\delta'\sigma) = (\sigma\delta\delta'\sigma)\delta'$ and 
(c) $\sigma\delta'\sigma\delta' = \delta'\sigma\delta'\sigma$ 
then $\sigma\delta\sigma\delta = \delta\sigma\delta\sigma$. \\
(iii) In the presentation of Theorem \ref{thm:presbg1}, we can 
replace relation $(CR3)$ by:
\\ $(CR3')$ $\delta_s \sigma_1 \delta_r \delta_{s} \sigma_1 =  \sigma_1 \delta_r\delta_{s} \sigma_1 \delta_s$ \\ 
for $1 \le r<s \le 2g+p-1$ with $(r,s)\neq (p+2i, p+2i+1)$, $0\leq i \leq g-1$.  \end{lem}   
\proof 
(i) Assume (a) $\delta\sigma\delta\delta'\sigma = \sigma\delta\delta'\sigma\delta$, 
(b) $\sigma\delta\sigma\delta = \delta\sigma\delta\sigma$ and (c) $\sigma\delta'\sigma\delta' = \delta'\sigma\delta'\sigma$. 
Then, $\delta'\sigma\delta\delta'\sigma = \delta^{-1}\sigma^{-1}\sigma\delta\delta'\sigma\delta\delta'\sigma = 
\delta^{-1}\sigma^{-1}\delta\sigma\delta\delta'\sigma\delta'\sigma = 
\delta^{-1}\sigma^{-1}\delta\sigma\delta\sigma\delta'\sigma\delta' = 
\delta^{-1}\sigma^{-1} \sigma\delta\sigma\delta \delta'\sigma\delta' = \sigma\delta \delta'\sigma\delta'$. \\
(ii) Assume (a) $\delta(\sigma\delta\delta'\sigma) = \sigma(\delta\delta'\sigma\delta)$, 
(b) $\delta'(\sigma\delta\delta'\sigma) = (\sigma\delta\delta'\sigma)\delta'$ and 
(c) $\sigma\delta'\sigma\delta' = \delta'\sigma\delta'\sigma$. 
Then \displaylines {\sigma\delta\sigma\delta = \sigma\delta\delta'\sigma\delta(\delta'\sigma\delta)^{-1}\sigma\delta = 
\delta\sigma\delta\delta'\sigma {\delta}^{-1}\sigma^{-1}{\delta'}^{-1}\sigma\delta =
\delta\sigma\delta(\sigma\delta'\sigma\delta'\sigma^{-1}{\delta'}^{-1}) {\delta}^{-1}\sigma^{-1}{\delta'}^{-1}\sigma\delta = 
\hfill\cr\hfill\delta\sigma\delta\sigma\delta'\sigma\delta'{\delta'}^{-1}\sigma^{-1}{\delta'}^{-1}{\delta}^{-1}\sigma^{-1}\sigma\delta = 
\delta\sigma\delta\sigma.}\noindent (iii) is a consequence of (i). \edim 

Since the relations of the presentation of $\bns$ are positive, one can define a monoid with the same presentation but as a monoid presentation. It is easy to see that the 
monoid we obtain doesnot inject in $\bns$, even if we add the relations of type $(CR3')$ to the presentation given in 
Theorem \ref{thm:presbg1}. In fact the following relations, $$(CR3)_{k} \ \ \delta_r \sigma_1 \delta_r \delta^k_{s} \sigma_1 = 
 \sigma_1 \delta_r\delta^k_{s} 
\sigma_1 \delta_r$$ for $1 \le r<s \le 2g+p-1$ with  $(r,s)\neq (p+2i, p+2i+1)$, $0\leq i \leq g-1$ and $k\in \mathbb{N}^*$,
 and $$(CR3')_k \ \ \delta_s \sigma_1 \delta^k_r \delta_{s} \sigma_1 =  \sigma_1 \delta^k_r\delta_{s} \sigma_1 \delta_s$$ 
for $1 \le r<s \le 2g+p-1$ with $(r,s)\neq (p+2i, p+2i+1)$, $0\leq i \leq g-1$ and $k\in \mathbb{N}^*$,  are true in $\bns$ 
for each positive integer $k$, but they are false in the monoid for $k$ greater than 1: no relation of the presentation can 
be applied to the left side of the equalities. Then starting from the left side of the equality for $k > 1$, we cannot obtain 
the right side of the equality by using the relations of the monoid presentation only.
\begin{quest} Let $B^*_n(\Sigma_{g,p})$ be the monoid defined by the presentation of Theorem \ref{thm:presbg1} with the extra 
relations $(CR3)_k$, $(CR3')_k$ for $k\in\mathbb{N}^*$. Is the canonical homomorphism $\varphi$ from $B^*_n(\Sigma_{g,p})$ 
to $\bns$ into ?     
\end{quest}
We remark that we can define a length function  $\ell$ on $B^*_n(\Sigma_{g,p})$: if $F^*$ is the free monoid based 
on $\sigma_1\cdots,\sigma_{n-1},\delta_1,\cdots, \delta_{2g}$, if  $l: F \to \mathbb{N}$ is the canonical length function 
and if $w\mapsto \overline{w}$ is the canonical morphism from $F^*$ onto $B^*_n(\Sigma_{g,p})$ then, for each $g$ 
in $B^*_n(\Sigma_{g,p})$, one has $sup\{l(w)\mid w\in F^*$ ; $\overline{w} = g\} < +\infty$; furthermore if 
we set $\ell(g) = sup\{l(w)\mid w\in F^*\}$, then for $g_1,g_2$ in $B^*_n(\Sigma_{g,p})$ 
we have $\ell(g_1g_2) \leq \ell(g_1)+\ell(g_2)$.\\       

Now, let us consider the particular case of planar surfaces.
\begin{prop}\label{propallcocksuite} Let $n,p$ be positive integers with $n\geq p-1$. Let $I\subset \{1,\cdots, n\}$
 with $Card(I) = p-1$.\\ Then $B_n(\Sigma_{0,p})$ admits the following presentation:\\ 
$\bullet$ Generators : $\sigma_1,\cdots,\sigma_{n-1}$ and $\rho_i$ for $i\in I$;\\$\bullet$ Relations:\\ \begin{tabular}{lll}
$(BR1)$ &$\ij  =  \ji$&for $1 \le i,j\le n-1$ with $| i-j | \ge 2$;\\
$(BR1)'$ &$\rho_r\rho_s = \rho_s\rho_r$&$r,s\in I$, $r\neq s$;\\
$(BR1)''$ &$\rho_r\sg_i = \sg_i\rho_r$&$r\in I$; $1 \le i\le n-1$, $i\ne r-1,r$.\\
$(BR2)$ &$\iji  = \jij$&for $1 \le i\le n-1$;\\
$(BR3)$ &$\sg_i\rho_r\sg_i\rho_r = \rho_r\sg_i\rho_r\sg_i$&$r\in I$, $i = r,r-1$;\\ 
$(BR3)'$ &$(\sg_{r-1}\sg_r)\rho_r\sg_{r-1}\rho_r = \rho_r(\sg_{r-1}\sg_r)\rho_r\sg_{r-1}$&$r\in I$ ; $r\neq 1,n$.
\end{tabular}\end{prop}
\proof Consider the presentation of Theorem \ref{thm:presbg1}. Let $I = \{r_1<r_2\cdots <r_{p-1}\}$ and 
set $\rho_{r_j} = (\sg_{r-1}\cdots \sg_1)\delta_{p-j}(\sg_{r-1}\cdots \sg_1)^{-1}$. Relations $(BR1)''$  
are equivalent to relations $(CR1)$ with $i\neq r$, by using braid relations. Using $(CR2)$, we get $(BR1)'\iff (CR3)$. 
Using braid relations $(BR1)$ and  $(BR2)$, the relations $(CR2)$ are equivalent to relations $(BR3)$ 
by conjugation by $(\sg_{r-1}\sg_r)\cdots(\sg_1\sg_2)$ and $(\sg_{r-2}\sg_{r-1})\cdots(\sg_1\sg_2)$ 
when $i = r-1$ and $i = r$ respectively . Now consider the relation $(CR1)$ for $i = r$. 
By conjugation by $\sg_{r-1}\cdots \sg_1$, we get $\sg_{r-1}\sg_r\sg_{r-1}^{-1}\rho_r = \rho_r\sg_{r-1}\sg_r\sg_{r-1}^{-1}$ and, 
then $\sg_{r-1}\sg_r\rho_r\sg_{r-1} = \sg_r\rho_r\sg_{r-1}\sg_r$. It follows that the relations of type $(CR1)$ for $i = r$  is 
equivalent to the relation $\sg_{r-1}\sg_r\rho_r\sg_{r-1} = \sg_r\rho_r\sg_{r-1}\sg_r$. This last relation is 
equivalent to relation $(BR3)'$ using relation $(BR3)$: \displaylines{\sg_{r-1}\sg_r\rho_r\sg_{r-1} = \sg_r\rho_r\sg_{r-1}\sg_r \iff 
\sg_{r-1}\sg_r\rho_r\sg_{r-1}\rho_r\sg_{r-1}  = \sg_r\rho_r\sg_{r-1}\sg_r\rho_r\sg_{r-1} \iff\hfill\cr\hfill
 \sg_{r-1}\sg_r\sg_{r-1}\rho_r\sg_{r-1}\rho_r = \sg_r\rho_r\sg_{r-1}\sg_r\rho_r\sg_{r-1} \iff
 \sg_r \sg_{r-1}\sg_r\rho_r\sg_{r-1} \rho_r= \sg_r\rho_r\sg_{r-1}\sg_r\rho_r\sg_{r-1} \iff\cr\hfill
 \sg_{r-1}\sg_r\rho_r\sg_{r-1}\rho_r= \rho_r\sg_{r-1}\sg_r\rho_r\sg_{r-1}.}
\edim
\begin{cor}(\cite{All} Table 1.1)\\ $B_n(\Sigma_{0,3})$ is isomorphic to the Affine Artin group of type $\tilde{B}(n+1)$ for $n\geq 2$.  
\end{cor}
\proof We apply Proposition \ref{propallcocksuite} with $I = \{1,n\}$. \edim

Recall that a monoid $M$ is cancellative if the property ``$\forall x,y,z,t \in M, (xyz = xtz)\Rightarrow (y=t)$'' holds in $M$.         
\begin{quest} Let $B^+_n(\Sigma_{0,p})$ the monoid defined by the presentation given in Proposition \ref{propallcocksuite}, 
considered as a monoid presentation.\\(i) Is the monoid $B^+_n(\Sigma_{0,p})$ cancellative ?\\(ii) is the natural homomorphism 
from $B^+_n(\Sigma_{0,p})$ to $B_n(\Sigma_{0,p})$  injective ?
\end{quest}
For $p = 1$ and $p=2$, the groups $B_n(\Sigma_{0,p})$ are isomorphic to the braid group $B_n$ and the Artin-Tits group of 
type $B$ respectively. Hence, the answer to above questions are positive.  In the case of $B_n(\Sigma_{0,3})$, 
the answers are also positive (see \cite{Cor} and \cite{Par}). Note that the relations of the presentation of $B_n(\Sigma_{0,p})$ 
 are homogeneous. Therefore we can define a length function $\ell$ on $B^+_n(\Sigma_{0,p})$ such that $\ell(g_1g_2) = \ell(g_1)+\ell(g_2)$ 
for every $g_1,g_2\in B^+_n(\Sigma_{0,p})$.   
\section{Braid groups on closed surfaces}
In this section, we consider braid groups on closed surfaces, that is without boundary components. In particular, we prove Corollaries \ref{thm:presbg2} 
and \ref{thm:presbg3}.
\begin{prop} \label{pressansbord}Let $n,g$ be positive integers. The group $B_n(\Sigma_{g,0})$ admits 
the following presentation:\\$\bullet$ Generators: $\sigma_1\cdots,\sigma_{n-1},\delta_1,\cdots,\delta_{2g}$;\\
$\bullet$ Relations\\ 
-Braid relations:
\noindent \begin{tabular}{lll}
$(BR1)$&$\sip \sjp = \sjp \sip$ &  for $| i-j | \ge 2$.\\
$(BR2)$&$\sip \siip \sip =  \siip \sip \siip$&$1\le i\le n-1$;
\end{tabular}\\
-Commutative relation between surface braids:\\
\begin{tabular}{lll}
$(CR1)$&$\sigma_i\delta_r = \delta_r\sigma_i$&$3\le i\le n-1; 1\le r\le 2g$;\\$(CR4)$&$\sigma_1\delta^2_r = 
\delta^2_r\sigma_1$&$1\le r\le 2g$;\\
&$\sigma_2\delta_{2r-1}\sigma_2 = \delta_{2r-1}\sigma_2\sigma_1$&$1\le r\le 2g$;\\&$\sigma_1\delta_{2r}\sigma_2=
 \sigma_2\delta_{2r}\sigma_1$&$1\le r\le 2g$;\end{tabular}\\
-Skew commutative relations on the handles:\\
(SCR2) $\sg_1\delta_r\delta_{r+2s}\sg_1 =\delta_{r+2s}\delta_{r}$ $1 \le r <r+2s\le 2g$;\\
$(SCR3)$ $\delta_{2r}\sg_1\delta_{2s-1}\delta_{2r}\sg_1 = \delta_{2s-1}\delta_{2r}^2$  $1\le s\le r\le g$;\\
$\sg_1\delta_{2s}\delta_{2r-1}\sg_1\delta_{2s} = \delta^2_{2s}\delta_{2r-1}$ $1\le s< r\le g$;\\
-Relation associated to the fundamental group of the surface:\\
(FGR) $(\sg_2\cdots\sg_{n-2}\sg_{n-1}^2\sg_{n-2}\cdots \sg_2)\sg_1\delta_1\delta_2\cdots\delta_{2g}\sg_1 = 
\delta_{2g}\cdots\delta_2\delta_1$.
\end{prop}
\proof
Starting from the presentation of Theorem \ref{deuxiemepresentation}, 
we set $\sigma_i = \theta_i^{-1}$, $\delta_{2r} = b_{2r}\theta_1^{-1}$ and $\delta_{2r-1} = \theta_1b^{-1}_{2r-1}$; 
we obtain easily the required presentation.
\edim

\begin{cor} \label{thm:presbg2}  Let $n$ and $g$ be positive integers with $g\geq 2$. Then, the group $B_n(\Sigma_{g,0})$ 
admits the following group presentation:\\
$\bullet$ Generators: $\genu,\delta_1,\cdots , \delta_{2g}$;\\
$\bullet$ Relations:\\
-Braid relations:\\
\begin{tabular}{lll}
$(BR1)$ &$\ij  =  \ji$&for $1 \le i,j\le n-1$ with $| i-j | \ge 2$;\\
$(BR2)$ &$\iji  = \jij$&for $1 \le i\le n-1$.\end{tabular}\\
- Commutative relations between surface braids:\\
\begin{tabular}{lll}
$(CR1)$&$\delta_r \sg_i= \sg_i \delta_r$  &for $2<i$; $1 \le r \le 2g$;\\
$(CR4)$&$\delta^2_r \sg_1= \sg_1 \delta^2_r$ &for $1 \le r \le 2g$;\\
&$\sg_2\delta_{2r-1}\sg_2 = \delta_{2r-1}\sg_2\sg_1$;\\
&$\sg_1\delta_{2r}\sg_2 = \sg_2\delta_{2r}\sg_1$;\\
$(CR5)$&$(\delta_{2r}\sg_1)(\delta_{2s-1}\delta_{2s}) = (\delta_{2s-1}\delta_{2s})(\delta_{2r}\sg_1) $&$1\le s<r\le g$;\\ 
&$(\sg_1\delta_{2r})(\delta_{2s}\delta_{2s-1}) = (\delta_{2s}\delta_{2s-1})(\sg_1\delta_{2r}) $&$1\le r<s\le g$;\\
&$(\delta_{2r}\sg_1)(\delta_{2r-1}\delta_{2s}) = (\delta_{2r-1}\delta_{2s})(\delta_{2r}\sg_1) $&$1\le s<r\le g$;\\
&$(\sg_1\delta_{2r-1})(\delta_{2s-1}\delta_{2r}) = (\delta_{2s-1}\delta_{2r})(\sg_1\delta_{2r-1}) $&$1\le r<s\le g$;
\end{tabular}\\
-Skew commutative relations on the handles\\ 
\begin{tabular}{lll}$(SCR2)$ &$\sg_1\delta_r\delta_{r+2s}\sg_1 =\delta_{r+2s}\delta_{r}$ &$1 \le r <r+2s\le 2g$;\end{tabular}\\
-Relation associated to the fundamental group of the surface\\
\begin{tabular}{lll}$(FGR)$ $(\sg_2\cdots\sg_{n-2}\sg_{n-1}^2\sg_{n-2}\cdots \sg_2)\sg_1\delta_1\delta_2\cdots\delta_{2g}\sg_1 = 
\delta_{2g}\cdots\delta_2\delta_1$.\end{tabular}
\end{cor}

\begin{cor} \label{thm:presbg3}  Let $n$ be positive integer. Then, 
the group $B_n(\Sigma_{1,0})$ admits the following group presentation:\\
$\bullet$ Generators: $\genu,\delta_1, \delta_{2}$;\\
$\bullet$ Relations:\\
-Braid relations:\\
\begin{tabular}{lll}
$(BR1)$ &$\ij  =  \ji$&for $| i-j | \ge 2$;\\
$(BR2)$ &$\iji  = \jij$&for $1 \le i\le n-1$.\end{tabular}\\
-Commutative relations between surface braids:\\
\begin{tabular}{lll}
$(CR1)$&$\delta_r \sg_i= \sg_i \delta_r$  &for $2<i$; $r = 1,2$;\\
$(CR4)$&$\delta^2_r \sg_1= \sg_1 \delta^2_r$ &$r = 1,2$;\\
&$\sg_2\delta_1\sg_2 = \delta_1\sg_2\sg_1$;\\
&$\sg_1\delta_2\sg_2 = \sg_2\delta_2\sg_1$;\end{tabular}\\
-Skew commutative relations on the handles:\\
$(SCR4)$ $\delta_2\sg_1\delta_1\delta_2\sg_1 = \delta_1\delta_2^2$;\\
-Relation associated to the fundamental group of the surface:\\
$(FGR)$ $(\sg_2\cdots\sg_{n-2}\sg_{n-1}^2\sg_{n-2}\cdots \sg_2)\sg_1\delta_1\delta_2\sg_1 = \delta_2\delta_1$.
\end{cor}
Corollary \ref{thm:presbg3} is a special case of Proposition \ref{pressansbord} when $g = 1$. When $g\geq 2$, 
Corollary \ref{thm:presbg2} follows from Proposition  \ref{pressansbord} using Lemma \ref{lemmebelow} below.

\begin{lem}\label{lemmebelow} (i) Let $G$ be a group and $\sigma,\delta_1,\delta_2,\delta_1',\delta_2'$ be in $G$ 
such that a) $\sigma\delta_i^2 = \delta_i^2\sigma$ for $i = 1,2$; b) $\sigma{\delta'}_i^2 = {\delta'}_i^2\sigma$ for
 $i = 1,2$; c) $g\delta'_2\delta_2g = \delta_2\delta'_2$ and d) $g\delta'_1\delta_1g = \delta_1\delta'_1$. 
Then,\\ (i) $\delta_2\sg\delta'_1\delta_2\sg = \delta'_1\delta_2^2 \iff \delta_2\sg \delta'_1\delta'_2= 
\delta'_1\delta'_2\delta_2\sg\iff \delta_1\delta_2\sg\delta'_1 = \sg\delta'_1\delta_1\delta_2$.\\(ii) $\delta_2\sg\delta_1\delta_2\sg = 
\delta_1\delta_2^2 \iff \delta_2\sg\delta_1\delta'_2 = \delta_1\delta'_2\delta_2\sg$.\\(iii) $\delta'_2\sg\delta'_1\delta'_2\sg = 
\delta'_1{\delta'}_2^2 \iff \delta_1\delta'_2\sg\delta'_1 = \sg\delta'_1\delta_1\delta'_2$. \end{lem}
\proof We prove under the hypothesis that 
$\delta_2\sg\delta'_1\delta_2\sg = \delta'_1\delta_2^2 \iff \delta_2\sg \delta'_1\delta'_2= \delta'_1\delta'_2\delta_2\sg$. 
The other cases are similar.\\ \displaylines{\delta_2\sg\delta'_1\delta_2\sg = \delta'_1\delta_2^2 \iff \delta_2\sg\delta'_1\delta_2 = 
\delta'_1\delta_2^2\sg^{-1}\iff \delta_2\sg\delta'_1\delta_2 = \delta'_1\sg^{-1}\delta_2^2\iff \delta_2\sg\delta'_1 = 
\delta'_1\sg^{-1}\delta_2\iff \hfill\cr\hfill \delta_2\sg\delta'_1\delta'_2 = 
\delta'_1\sg^{-1}\delta_2\delta'_2\iff \delta_2\sg \delta'_1\delta'_2= \delta'_1\delta'_2\delta_2\sg.}
\edim

\begin{lem} Let $n,g$ be positive integers. Consider the group $B_n(\Sigma_{g,0})$ and the presentation of Proposition \ref{pressansbord}.
 Then, for every $1\le r,s\le g$, $$\delta_{2s-1}\delta^2_{2r}\delta_{2s-1}  = \delta_{2r}\delta^2_{2s-1}\delta_{2r}.$$
\end{lem} 
\proof Let $1\le s\le r\le g$. From  the relations of type $(SCR3)$ and the relations of the first type of $(CR4)$, 
it follows $\delta_{2s-1}\delta_{2r}\sg_1\delta_{2s-1}\delta_{2r}\sg_1 = \delta_{2s-1}^2\delta_{2r}^2 = 
\sg_1\delta_{2s-1}\delta_{2r}\sg_1\delta_{2s-1}\delta_{2r}$. Hence, $\delta_{2s-1}^2\delta_{2r}\delta_{2s-1}^{-1} = 
\sg_1\delta_{2s-1} \delta_{2r} \sg_1 = \delta_{2r}^{-1}\delta_{2s-1}\delta_{2r}^2$ and then $\delta_{2s-1}\delta^2_{2r}\delta_{2s-1}  =
 \delta_{2r}\delta^2_{2s-1}\delta_{2r}$.\\If $1\le r < s \le g$ then we proceed in the same way, using that 
$\sg_1\delta_{2s}\delta_{2r-1}\sg_1\delta_{2s} = \delta^2_{2s}\delta_{2r-1}$.  
\edim
\begin{quest} Consider the monoids defined by the presentation given in Theorem 3.1, Corollary 3.2  or in Corollary 3.3. 
Are they cancellative ? do they embed in $B_n(\Sigma_{g,0})$ ? 
\end{quest}

\section{Braid group on two strands on the torus}
\subsection{The word problem and the conjugacy problem}
In this section we solve the word problem and the conjugacy problem for the special case when  $g = 1$ and $n = 2$ by using a 
presentation derived from the one obtained in the previous section.\\ As a consequence of Corollary \ref{thm:presbg3} we have:
\begin{cor} \label{proppresb210} $B_2(\Sigma_{1,0})$ admits the group presentation:
$$B_2(\Sigma_{1,0} = \langle \sigma_1,\delta_1, \delta_{2}\mid 
 \delta^2_1 \sg_1= \sg_1 \delta^2_1; \delta^2_2 \sg_1= \sg_1 \delta^2_2; \delta_2\sg_1\delta_1\delta_2\sg_1 = 
\delta_1\delta_2^2; \sg_1\delta_1\delta_2\sg_1 = \delta_2\delta_1\rangle.$$
\end{cor}
\begin{lem} \label{lempresb210} Let $G$ be group and $\sigma_1,\delta_1, \delta_{2}$ be in $G$ such that  
a) $\delta^2_1 \sg_1= \sg_1 \delta^2_1$~; b) $\delta^2_2 \sg_1= \sg_1 \delta^2_2$~; 
c) $\sg_1\delta_1\delta_2\sg_1 = \delta_2\delta_1$. 
Then $\delta_2\sg_1\delta_1\delta_2\sg_1 = \delta_1\delta_2^2 \iff (\delta^2_1 \delta_2= \delta_2 \delta^2_1$ 
and $\delta^2_2 \delta_1 = \delta_1 \delta^2_2)$.
\end{lem}
\begin{cor}\label{lemmurKur}(\cite{MurKur}, Chapter 11, Exercises 5.2 and 6.3) The group $B_2(\Sigma_{1,0})$ 
admits the two following group presentations: $$ B_2(\Sigma_{1,0}) = \langle \sigma_1,\delta_1, \delta_{2}\mid \delta^2_1 \sg_1=
 \sg_1 \delta^2_1~; \delta^2_2 \sg_1= \sg_1 \delta^2_2~; \delta^2_1 \delta_2= \delta_2 \delta^2_1~; \delta^2_2 \delta_1 = 
\delta_1 \delta^2_2~; \sg_1\delta_1\delta_2\sg_1 = \delta_2\delta_1 \rangle.$$ \begin{eqnarray}\label{preb2} B_2(\Sigma_{1,0}) = 
\langle a,b,c\mid a^2b= ba^2; b^2a = ab^2; a^2c= ca^2; b^2c = cb^2;a^2b^2 = c^2 \rangle.\end{eqnarray}\end{cor}
\proof (i) follows from Corollary \ref{proppresb210} and Lemma \ref{lempresb210}. For ii), we set $a = \delta_2$, $b = \delta_1$ and
 $c = \delta_2\delta_1\sg_1^{-1}$ as suggested in \cite{MurKur} Chapter 11, Exercise 6.3. \edim

Using the presentation (\ref{preb2}), we are able to solve the word problem and the conjugacy problem 
in $B_2(\Sigma_{1,0})$. Considering (\ref{preb2}), for $x = a$, or $x = b$ we can define a weight homomorphism 
of groups $\ell_{\hat{x}} : B_2(\Sigma_{1,0}) \to \mathbb{Z}$ such that $\ell_{\hat{x}}(x) = 0$ 
and $\ell_{\hat{x}}(y) = 1$ for $y\in\{a,b,c\}$ and $y\neq x$.\\ In the following we denote by $F(a,b,c)$ the free group 
based on $\{a,b,c\}$. We denote by $W(a,b,c)$ the Coxeter group associated to $F(a,b,c)$ and defined by 
$W(a,b,c) = \langle a,b,c\mid a^2 = b^2 = c^2 = 1\rangle$. If $w$ is in $F(a,b,c)$ we denote by $\overline{w}$ 
its image in $B_2(\Sigma_{1,0})$. Considering (\ref{preb2}), there exists a morphism 
$p : B_2(\Sigma_{1,0})\to W(a,b,c)$ that sends $x\in \{a,b,c\}$ on $x$. Note that the canonical 
morphism from $F(a,b,c)$ onto $W(a,b,c)$ factorises through $p$.\\  We denote by $L_{a,b}$ 
the set-map from $B_2(\Sigma_{1,0})$ to  $F(a,b,c)$  defined by $L_{a,b}(g) = a^{\ell_{\hat{b}}(g)}b^{\ell_{\hat{a}}(g)}$ 
for $g$ in $B_2(\Sigma_{1,0})$.  If $w$ is in $F(a,b,c)$, we write, by abuse of notation, $p(w)$ and $L_{a,b}(w)$ 
for $p(\overline{w})$ and $L_{a,b}(\overline{w})$ respectively. 
\begin{prop}
(i) The center $Z(B_2(\Sigma_{1,0}))$ of the group $B_2(\Sigma_{1,0})$ is a free Abelian group based on $a^2$ and $b^2$. 
Furthermore, for each element $g$ of $Z(B_2(\Sigma_{1,0}))$, the word $L_{a,b}(g)$ is a representing element of $g$.\\
(ii) The group $B_2(\Sigma_{1,0})$ is a central extension of $W(a,b,c)$.\\ 
In other words, the sequence $1\to Z(B_2(\Sigma_{1,0})) \to B_2(\Sigma_{1,0}) \to W(a,b,c)\to 1$ is exact. 
\end{prop}
\proof We remark that the presentation (\ref{preb2}) implies that $a^2$, $b^2$  and $c^2$ 
are in  $Z(B_2(\Sigma_{1,0}))$ and that $W(a,b,c) = B_2(\Sigma_{1,0})/a^2 = b^2 = 1$. 
Since the center of $W(a,b,c)$ is trivial, (ii) follows. As a consequence, 
we get that $a^2$ and $b^2$ generated the Abelian group $Z(B_2(\Sigma_{1,0}))$. Now, let $g$ belong 
to  $Z(B_2(\Sigma_{1,0}))$. Each word $w$  that represents $g$ and written on the letters $a^2$, $b^2$, 
and their inverses, can be modify in order to obtain $L_{a,b}(g)$ by using the relations $a^2b^2 = b^2a^2$ 
and the relations of $a^2$ and $b^2$ with their respective inverses. Hence,  $Z(B_2(\Sigma_{1,0}))$ is a 
free Abelian group based on $a^2$ and $b^2$.  

\edim
\begin{cor}\label{propfinalpourcasparticuliergroup}  
Let $w$  be in $F(a,b,c)$; then  $\overline{w} = 1\iff (\  L_{a,b}(w) = 1$ and $p(w) = 1\ )$.\end{cor}

\proof Assume  $L_{a,b}(w) = 1$ and $p(w) = 1$. Since $p(\overline{w}) = 1$, the element 
$\overline{w}$ is in the center of $B_2(\Sigma_{1,0})$. But in that case, $L_{a,b}(w) = 1$ 
represents $\overline{w}$. Then $\overline{w} = 1$. 
\edim
\begin{cor} The word problem in $B_2(\Sigma_{1,0})$ is solvable.\end{cor}
\proof The word problem is solvable in the free group $F(a,b,c)$ and in the Coxeter group $W(a,b,c)$. 
Then the claim follows from Corollary \ref{propfinalpourcasparticuliergroup}.\edim

\begin{cor}\label{lempourconj} Let $g,h$ be  in $F(a,b,c)$; then, $\overline{g} = \overline{h} \iff p(g) = p(h)$ 
and $L_{a,b}(g) = L_{a,b}(h)$. 
\end{cor}

\proof  Assume $p(g) = p(h)$  and $L_{a,b}(g) = L_{a,b}(h)$.
Then the element $\overline{g}\overline{h}^{-1}$ is in the center of $B_2(\Sigma_{1,0})$.
Since $L_{a,b}(g) = L_{a,b}(h)$ it follows that $\ell_{\hat{a}}(g)= \ell_{\hat{a}}(h)$
and $\ell_{\hat{b}}(g)= \ell_{\hat{b}}(h)$. Therefore $L_{a,b}(gh^{-1})=1$ and thus  $\overline{g}=\overline{h}$.
\edim

We denote by $F^+(a,b,c)$ the free monoid based on $\{a,b,c\}$. It is a submonoid of $F(a,b,c)$. 
If $w$ is in $W(a,b,c)$, there exists a unique element $[w]$ in $F^+(a,b,c)$ of minimal length such that 
its image in $W(a,b,c)$ is $w$.  By construction, for each $w$ in $W(a,b,c)$ the element $p(\overline{[w]})$ is equal to $w$. 
As a consequence, the map sending each $w$ in $W(a,b,c)$ on $\overline{[w]}$ is injective. For short, 
we will write $[w]$ for $\overline{[w]}$. 
\begin{cor} (i) Let $g,h$ be in $B_2(\Sigma_{1,0})$; 
then,$$(\exists r\in B_2(\Sigma_{1,0}), rgr^{-1} = h )\iff ( L_{a,b}(g) = L_{a,b}(h)\textrm{ and } \exists w\in W(a,b,c), wp(g)w^{-1} = 
p(h)).$$ Furthermore, if the right side  holds, then $[w]g[w]^{-1} = h$.   
\end{cor}
\proof
 The side ``$\Rightarrow$'' is clear with $w = p(r)$. Assume conversely that $ L_{a,b}(g) = L(a,b)(h)$  and  
$\exists w\in W(a,b,c),\ wp(g)w^{-1} = p(h)$. Since   $wp(g)w^{-1} = p(h)$, we have $p([w]g[w]^{-1}) = p(h)$. 
But $L_{a,b}([w]g[w]^{-1}) = L_{a,b}(g) = L_{a,b}(h)$ and $[w]g[w]^{-1}$ = $h$ by Corollary \ref{lempourconj}.
\edim
\begin{cor}The conjugacy problem in $B_2(\Sigma_{1,0})$ is solvable.\end{cor}
\proof The conjugacy problem is solvable in each Coxeter group (see \cite{Kra}).\edim
\subsection{The Garside method and complete presentation}
In order to solve the word problem and the conjugacy problem in $B_2(\Sigma_{1,0})$, we can try to use the 
method used by Garside to solve the word problem and the conjugacy problem, that is to find a Garside structure for $B_2(\Sigma_{1,0})$.
Let us remark that surface braid groups on surfaces of genus greater than $1$ have trivial center (see \cite{ParRol}) and then
they cannot be Garside groups. 
Recall that in  a monoid $M$ we say that $a$ left-divides $b$ if $b = ac$ for some $c$ in $M$. We say, in a similar way, 
that $a$ right-divides $b$ when $b = ca$ for some $c$ in $M$. An element $\Delta$ of $M$ is said to be balanced 
when its set of left-divisors is equal to its set of right-divisors.  We denote by  $B^+_2(\Sigma_{1,0})$ 
 the monoid defined by the presentation (\ref{preb2}), but considered as a monoid presentation. 
Then in  $B^+_2(\Sigma_{1,0})$ the element $c^2$ is balanced. Furthermore its set $D(c^2)$ of divisors generates $B^+_2(\Sigma_{1,0})$. 
Nevertheless, $B^+_2(\Sigma_{1,0})$ fails to be a Garside monoid with $c^2$ for Garside element (see \cite{Deh} for a definition) 
because it is not a lattice for left-divisibility: $a$ and $b$ have two  distinct minimal common multiples, namely $ab^2$ and $ba^2$.  
Anyway, as shown in \cite{Deh} Section 8, part of the results established for Garside groups 
still hold, as we will see in Lemma \ref{corpropfinalpouras}.

Let $\iota: B_2^+(\Sigma_{1,0})\to  B_2(\Sigma_{1,0})$ be the canonical homomorphism of monoids. 
By abuse of notation, we denote by  $\ell_{\hat{a}}$ and  $\ell_{\hat{b}}$  the morphisms $\ell_{\hat{a}}\circ \iota$ 
and  $\ell_{\hat{b}}\circ \iota$ respectively. We remark that $z\mapsto\overline{z}$ factorises through $\iota$. 
By abuse of notation, we denote by  $z\mapsto\overline{z}$ and by $w\mapsto\overline{[w]}$ the factorizations. 
Then we have  $\overline{w} =\iota(\overline{w})$. As before, we write $[w]$ for $\overline{[w]}\in  B_2^+(\Sigma_{1,0})$. 
We remark that Corollary 4.5 and 4.6 still hold if we consider $\overline{w}$ in $B_2^+(\Sigma_{1,0})$ 
and $g,h$ in  $B_2^+(\Sigma_{1,0})$. As a consequence, we have the following result: 
\begin{lem} \label{corpropfinalpouras}(i) $ B^+_2(\Sigma_{1,0})$ is cancellative and the canonical 
morphism $\iota : B^+_2(\Sigma_{1,0}) \to  B_2(\Sigma_{1,0})$ is into. \\
(ii) $\forall G \in  B_2(\Sigma_{1,0})$, $\exists ! j\in\mathbb{Z}, \exists ! g\in  B^+_2(\Sigma_{1,0})$  
such that $G = c^{-2j} \iota(g)$ and $c^2$ doesnot divide $g$.
\end{lem} 
\proof (i) Let $h,g_1,g_2$ be in  $B^+_2(\Sigma_{1,0})$ such that $hg_1 = hg_2$. 
Then $\ell_{\hat{a}}(hg_1) = \ell_{\hat{a}}(h)+ \ell_{\hat{a}}(g_1)$ 
and $\ell_{\hat{a}}(hg_2) = \ell_{\hat{a}}(h)+ \ell_{\hat{a}}(g_2)$. 
Hence we have $\ell_{\hat{a}}(g_1) = \ell_{\hat{a}}(g_2)$. In the same way we get  $\ell_{\hat{b}}(g_1) = \ell_{\hat{b}}(g_2)$, 
and also $p(g_1) = p(g_2)$ in the group $W(a,b,c)$. Then  $g_1 = g_2$. We proceed in the say way if $g_1h = g_2h$. 
The other results are consequences of the Garside like structure as proved in Proposition 8.10 of \cite{Deh}. \edim

In the following, we identify $ B^+_2(\Sigma_{1,0})$ with its image in $ B_2(\Sigma_{1,0})$. In order to solve 
the word problem in $B_2(\Sigma_{1,0})$, it is then enough to solve the words problem in $B^+_2(\Sigma_{1,0})$. Then, 
using the following proposition, we obtain another solution to the word problem.
 
\begin{prop}\label{propfinalpourcasparticulier}  Let $g$ be in $B^+_2(\Sigma_{1,0})$; then, there exist a 
unique pair $(h,l)$ in $\mathbb{N}^2$, and a unique $w$ in $W(a,b,c)$ such that $g = a^{2k}b^{2l}[w]$.
\end{prop}
\proof  Since $c^2 = a^2b^2$ and that $a^2$, $b^2$ are in the center of $Z(B^+_2(\Sigma_{1,0}))$, it 
follows that we can write $g = a^{2k}b^{2l}[w]$ with $k,l$ in $\mathbb{N}$ and  $w$ in $W(a,b,c)$. 
Assume that  $g = a^{2i}b^{2j}[z]$ for some $i,j$ in $\mathbb{N}$ and  $z$ in $W(a,b,c)$. We have (1) $ w = p(g)  z$ and 
then $[w] = [z]$; (2) $k = \frac{\ell_{\hat{b}}(g) - \ell_{\hat{b}}([w])}{2}= i$; 
(3) $l = \frac{\ell_{\hat{a}}(g) - \ell_{\hat{a}}([w])}{2} = j$. Then the decomposition of $g = a^{2k}b^{2l}[w]$ is unique.\edim

If we want to solve the conjugacy problem by using the idea of Garside, we need to understand the normal form of $B^+_2(\Sigma_{1,0})$ 
as defined in Definition 7.2 of \cite{Deh}. This lead us to the notion of complete presentation as defined in \cite{Deh2}. 
Let $S$ be a finite set, and $S^*$ be the free monoid based on $S$. We denote by $\epsilon$ the empty word. Let $B$ be a 
monoid with presentation $(S,\mathcal{R})$. We write $w\equiv w'$ if the word $w,w'$ of $S^*$ have the same image in $B$. 
Let $w,w'$ be two words in $(S\cup S^{-1})^*$, where $S^{-1}$ is a disjoint copy of $S$. We say that $w$ reverses in $w'$, 
and write $w\curvearrowright w'$ if $w'$ is obtained from $w$ by a finite sequence of the following steps:  
deleting some $u^{-1}u$ for some $u\in S$ or replacing some subword $u^{-1}v$ where $u,v$ are in $S$, 
with a a word $v'{u'}^{-1}$ such  that $uv' = vu'$ is a relation of $\mathcal{R}$.   

\begin{df}[\cite{Deh2} Definition 2.1 and Proposition 3.3] Let $B$ be a monoid with presentation $(S,\mathcal{R})$. 
We say that the presentation $(S,\mathcal{R})$ is complete if 
$$\forall u,v\in S^*,\  (u\equiv v \iff u^{-1}v\curvearrowright \epsilon).$$  
\end{df}

For instance the classical presentation of each Artin-Tits monoid is complete. The definition of 
complete presentation is easy to understand. Nevertheless, it is not  easy to verify that a given presentation 
is complete. In \cite{Deh2}, Dehornoy gives a semi-algorithmical method in order to decide if a given presentation is complete. 
Semi-algorithmical means that when the process finishes, it gives an answer, but it is possible that it doesnot finish.
 We do not explain this technical method, named {\it the cube condition}, but refer to Definition 3.1 and Figure 3.1 of \cite{Deh2}.    
  
Applying the cube condition process, it is quiet clear that the presentation (\ref{preb2}) of the monoid $B_2(\Sigma_{0,1})$ 
is not complete and that we must add to the presentation the relation ``$b^2a^2 = c^2$ '' if we want to expect that the 
presentation is complete.
\begin{quest} Is the presentation \begin{eqnarray}\langle a,b,c\mid a^2b= ba^2; b^2a = ab^2;
\label{presentcomplete} a^2c= ca^2; b^2c = cb^2 ; a^2b^2 = c^2 ;  b^2a^2 = c^2 \rangle^+ \end{eqnarray} a 
complete presentation of the monoid $B_2(\Sigma_{1,0})$ ? 
In other words, does this presentation verify the cube condition ?
\end{quest} 
\begin{quest} Is the monoid  presentation $\langle a,b\mid ab^2 = b^2a ; ba^2 = a^2b\rangle^+$ complete ? In other words, 
does this presentation verify the cube condition ?
\end{quest}
A positive answer to Question 5 seems to be crucial in order to state the interest of the method of completeness. 
\appendix
\section{Presentations of surface braid groups}
\begin{thm}[\cite{Bel} \label{thm:pres} Theorem 1.1] Let $n$, $p$ be  positive integers and $g$  a non negative integer. 
Let $g$ be a non negative integer.
The group $\bns$  admits the following group presentation:

\noindent $\bullet$ Generators: $\gend$.

\noindent $\bullet$  Relations:
\begin{itemize}
\item  braid relations:

\noindent \begin{tabular}{lll}
$(BR1)$&$\sip \sjp = \sjp \sip$ &  for $| i-j | \ge 2$.\\
$(BR2)$&$\sip \siip \sip =  \siip \sip \siip$&$1\le i\le n-1$;
\end{tabular}\\
\item  mixed  relations:

\noindent\begin{tabular}{lll}
(R1) &$a_r \sg_i=\sg_i  a_r$&for $1 \le r \le g;\; i\not= 1$;\\
&$b_r \sg_i=\sg_i  b_r$ &for $1 \le r \le g$; $i\not= 1$; \\
(R2)&$\sigma_1^{-1} a_r \sigma_1^{-1} a_{r}= a_r \sigma_1^{-1} a_{r}\sigma_1^{-1}$ &for $1 \le r \le g$; \\
&$\sigma_1^{-1} b_r \sigma_1^{-1} b_{r}= b_r \sigma_1^{-1} b_{r}\sigma_1^{-1}$&for $1 \le r \le g$;\\
(R3)&$\sigma_1^{-1} a_{s} \sigma_1 a_r = a_r \sigma_1^{-1} a_{s} \sigma_1$&for   $s < r$; \\
&$\sigma_1^{-1} b_{s} \sigma_1 b_r = b_r \sigma_1^{-1} b_{s} \sigma_1$&for $s < r$;\\
&$\sigma_1^{-1} a_{s} \sigma_1 b_r = b_r \sigma_1^{-1} a_{s} \sigma_1$&for $s < r$ \\
&$\sigma_1^{-1} b_{s} \sigma_1 a_r = a_r \sigma_1^{-1} b_{s} \sigma_1$&for $s < r$\\
(R4)&$\sigma_1^{-1} a_r \sigma_1^{-1} b_{r} = b_{r} \sigma_1^{-1} a_r \sigma_1$ &for $1 \le r \le g$;\\ 
(R5)&$z_j\sg_i= \sg_i z_j$  &for  $i \not=1, j=1, \dots, p-1$;\\ 
(R6)&$\sigma_1^{-1} z_i \sigma_1 a_r =a_r\sigma_1^{-1} z_i\sigma_1$&for $1 \le r \le g;\; i=1,
\dots, p-1$; \\
&$\sigma_1^{-1} z_i \sigma_1 b_r =b_r\sigma_1^{-1} z_i\sigma_1$  &for $1 \le r \le g;\; i=1,
\dots, p-1$\\
(R7)&$\sg_{1}^{-1} z_j \sg_{1} z_l=z_l\sg_{1}^{-1} z_j \sg_{1}$ &for $j=1, \dots, p-1 $, $j<l$;\\
(R8)&$\sg_{1}^{-1} z_j \sg_{1}^{-1} z_j=z_j\sg_{1}^{-1}  z_j \sg_{1}^{-1}$&for $ j=1, \dots, p-1$.
\end{tabular}
\end{itemize}
\end{thm}
\begin{thm}[\cite{Bel} Theorem A.4]\label{deuxiemepresentation} Let $\Sigma_{0,g}$ be a closed surface of positive genus $g$. 
Then $B_n(\Sigma_{0,g})$ admits the following presentation:\\$\bullet$ Generators: 
$\theta_1,\cdots \theta_{n-1}$ $b_1\cdots b_{2g}$;\\$\bullet$ Relations:\\ 
\noindent \begin{tabular}{lll}
$(BR1)$&$\theta_i\theta_j=\theta_j\theta_i$ &for $2\le | i-j |$.\\
$(BR1)$&$\theta_i \theta_{i+1} \theta_i= \theta_{i+1} \theta_i \theta_{i+1}$&for $1\le i\le n-2$.\end{tabular}\\
\noindent \begin{tabular}{lll}
$(R1)$&$b_r\theta_i = \theta_i b_r$&$1\le r \le 2g$; $i\neq 1$;\\
$(R2)$&$b_s\theta_1^{-1}b_r\theta_1^{-1} = \theta_1 b_r\theta_1^{-1}b_s$&$1\le s<r \le 2g$;\\
$(R3)$&$b_r\theta_1^{-1}b_r\theta_1^{-1} = \theta_1^{-1}b_r\theta_1^{-1}b_r$&$1\le r\le 2g$;
\end{tabular}\\
$(TR)$ $b_1b_2^{-1}\cdots b_{2g-1}b_{2g}^{-1}b_1^{-1}b_2\cdots b_{2g-1}^{-1}b_{2g} = 
\theta_1\cdots \theta_{n-2}\theta^2_{n-1}\theta_{n-2}\cdots\theta_1$.
\end{thm}

\end{document}